\documentclass[letterpaper, 10 pt, conference]{ieeeconf}

\IEEEoverridecommandlockouts                              
 \pdfoutput=1
\usepackage{amsmath, amssymb}
\usepackage{amsfonts}
\usepackage{graphicx}
\usepackage{enumerate}
\usepackage{ifthen}
\usepackage{multicol}
\usepackage{float}
\usepackage{cite}
\usepackage{fancyhdr}
\usepackage[usenames, dvipsnames]{color}
\graphicspath{{figures/}}
\usepackage{mathtools}
\usepackage{multicol}
\usepackage{xcolor}

 \usepackage{cancel} 
 \usepackage[normalem]{ulem}

\newtheorem{theorem}{Theorem}
\newtheorem{lemma}[theorem]{Lemma}
\newtheorem{proposition}[theorem]{Proposition}

\newtheorem{rem}{Remark}
\newtheorem{example}{Example}

\newtheorem{definition}{Definition}

\newboolean{showcomments}
\setboolean{showcomments}{true}

\newcommand{\todo}[1]{  \ifthenelse{\boolean{showcomments}}
{\textcolor{ForestGreen}{TO DO:  #1}}{}}
\newcommand{\suggest}[1]{\ifthenelse{\boolean{showcomments}}
{\textcolor{Orange}{(Suggestion: #1)}}{}}
\newcommand{\alain}[1]{\ifthenelse{\boolean{showcomments}}
{\textcolor{Blue}{(Alain says: #1)}}{}}
\newcommand{\jonas}[1]{\ifthenelse{\boolean{showcomments}}
{\textcolor{ForestGreen}{(Jonas says: #1)}}{}}
\newcommand{\kristian}[1]{\ifthenelse{\boolean{showcomments}}
{\textcolor{Blue}{(Kristian says: #1)}}{}}
\newcommand{\emma}[1]{\ifthenelse{\boolean{showcomments}}
{\textcolor{VioletRed}{(Emma says: #1)}}{}}
\newcommand{\ifneeded}[1]{\ifthenelse{\boolean{showcomments}}
{\textcolor{Gray}{#1}}{}}

\newboolean{showedit}
\setboolean{showedit}{true}
\newcommand{\edit}[1]{\ifthenelse{\boolean{showedit}}
{\textcolor{Blue}{#1}}{}}
\newcommand{\draft}[1]{\ifthenelse{\boolean{showedit}}
{\textcolor{gray}{#1}}{}}

\usepackage{array}
\newcolumntype{L}[1]{>{\raggedright\let\newline\\\arraybackslash\hspace{0pt}}m{#1}}
\newcolumntype{C}[1]{>{\centering\let\newline\\\arraybackslash\hspace{0pt}}m{#1}}
\newcolumntype{R}[1]{>{\raggedleft\let\newline\\\arraybackslash\hspace{0pt}}m{#1}}

\usepackage{physics} 

\usepackage[font = small]{caption}
\usepackage{subcaption}

\usepackage{tikz}
 \usetikzlibrary{plotmarks}

\usepackage{tikz}
\usepackage{pgfplots}
\pgfplotsset{compat=newest}
\usetikzlibrary{patterns}
\usetikzlibrary{decorations.text}
\usepgfplotslibrary{fillbetween}

\renewcommand{\baselinestretch}{0.985}

\providecommand{\figref}{}
\renewcommand{\figref}[1]{Fig.~\ref{fig:#1}}
\providecommand{\secref}{}
\renewcommand{\secref}[1]{Sec.~\ref{sec:#1}}

\usepackage{lipsum}
\usepackage{mathtools}

\title{\LARGE \bf Input-Output Pseudospectral Bounds for Transient Analysis of\\ Networked and High-Order Systems}

\author{ {Jonas Hansson and Emma Tegling} 
 \thanks{The authors are with the Department of Automatic Control,
        Lund University, Lund, Sweden. Email: \{{\tt\small{jonas.hansson, emma.tegling}\}@control.lth.se}}\thanks{This work was partially funded by Wallenberg AI, Autonomous Systems and Software Program (WASP) funded by the Knut and Alice Wallenberg Foundation and the Swedish Research Council through Grant 2019-00691. }
}

\begin{document}
\maketitle
\begin{abstract}
   Motivated by a need to characterize transient behaviors in large network systems in terms of relevant signal norms and worst-case input scenarios, we propose a novel approach based on existing theory for matrix pseudospectra. We extend pseudospectral theorems, pertaining to matrix exponentials, to an input-output setting, where matrix exponentials are pre- and post-multiplied by input and output matrices. Analyzing the resulting transfer functions in the complex plane allows us to state new upper and lower bounds on system transients. These are useful for higher-order matrix differential equations, and specifically control of double-integrator networks such as vehicle formation problems. Therefore, we illustrate the theory's applicability to the problem of vehicle platooning and the question of string stability, and show how unfavorable transient behaviors can be discerned and quantified directly from the input-output pseudospectra. 
\end{abstract}

\section{Introduction}
Characterizing dynamic properties of systems with structure, in particular, network structure, is a long-standing problem in the field. While questions of stability and convergence have dominated the literature since the early works~\cite{FaxMurray,OlfatiSaber2004}, important questions pertaining to the performance and robustness of network systems are increasingly gaining attention. For example, \cite{Bamieh2012} and later \cite{SiamiMotee2015, Tegling2019TAC} have described fundamental limitations to the performance of large networks subject to structural (sparsity) constraints, stated in terms of system norms. 

A particular area where dynamic behaviors have received more attention is that of vehicle platooning, that is, the control of strings of vehicles, see~\cite{LevineAthans1966,chu1974decentralized_control} for early works. Here, it is fundamentally important to prevent disturbance propagation through the string (to avoid collisions!), and therefore, to have uniform bounds on error amplifications during transients. 
This has motivated the notion of \textit{string stability}, see e.g.,~\cite{swaroop1996stringstability,seiler2004disturbancep_propagation} or~\cite{studli2017StringConcepts, FENG2019StringDefs} for more recent surveys. Conditions for string stability fall, roughly speaking, into two categories: 1) bounding the amplification of a disturbance from vehicle $i$ to vehicle $j$, or 2) requiring that bounded initial errors lead to bounded output errors, independently of the string length. The choice of signal norms, however, is central for the bounds in this literature, and the interpretations they allow for.  Many works have done analyses based on $ \mathcal{L}_2 $ to  $ \mathcal{L}_2 $ string stability, see \cite{HERMAN2017diffasym,seiler2004disturbancep_propagation, ploeg2014Lpstability} while the, as argued e.g. in~ \cite{feintuch2012}, possibly more important $\mathcal{L}_\infty$ to $ \mathcal{L}_\infty$ disturbance amplification has received significantly less attention even if considered in~\cite{swaroop1996stringstability,chu1974decentralized_control}. In this work, we shed light on a new approach to analyzing such bounds for input-output systems in general, and networks and vehicle strings in particular. 

This approach takes off from the literature on \textit{pseudospectra}. Pseudospectra, which complement spectral analysis of linear systems, especially for those with non-normal operators, have seen usage in describing the transient behavior of both differential and difference equations. The works are too numerous to mention, but we refer to~\cite{trefethenbook} for an excellent textbook on the subject. Through pseudospectra one can state lower and upper bounds on the transient of the exponential matrix, i.e., on $\sup_{t\geq0} \|e^{t\mathcal{A}}\|$, and thereby on the solution to a linear differential equation. In other words, on the transient response of the internal states of a linear system. The most famous such bounds are given by the Kreiss theorem~\cite{Kreiss1962}. However, in control, and in particular, network applications including vehicle platooning, we are not necessarily interested in the transients of the internal states. For instance, vehicular formation dynamics tend to have a double integrator rendering certain internal states unbounded, while inter-vehicular distances may be well-behaved. To cope with this one can incorporate measurement and input matrices $\mathcal{C},\mathcal{B}$ and then bound $\sup_{t\geq0} \|\mathcal{C}e^{t\mathcal{A}}\mathcal{B}\|$ instead. 

The extension of pseudospectral bounds to such an input-output setting is the main focus of the present work. For this purpose we will define a notion of \textit{input-output pseudospectra}. These
will, in the case of higher-order systems (by which we mean systems with more than one integrator), become closely related to \textit{structured} pseudospectra, which have been studied in~\cite{Structured2001Tisseur,lancaster2005matrix_polynomials} and applied to mechanical systems in~\cite{GREEN2006transient_response}. In these works the main focus has been on the robustness of solutions to matrix polynomial equations including the quadratic eigenvalue problem. The related analysis of transient behavior of $\|\mathcal{C}e^{tA}\mathcal{B}\|$ has, to the best of our knowledge, barely received attention, though some structured Kreiss-like theorems were proven in~\cite{matsuo1994maxamplitude,plischke2005transient}. 

This paper aims to highlight the potential usefulness of the pseudospectral framework for networked systems and systems with higher-order dynamics. Platooning, where vehicles are modeled as double integrators (the acceleration is actuated), and which have a string network topology, is a prototypical example. We first generalize certain key results from \cite{trefethenbook} to an input-output setting. Furthermore, we use complex analysis to derive new upper bounds on the 
transients of state space realizations, which are especially useful for systems that have high-order dynamics. The generalizations lead to lower and upper bounds on the transient $\sup_{t\geq0} \|\mathcal{C}e^{t\mathcal{A}}\mathcal{B}\|$, which under given input scenarios imply bounds on the output $\sup_{t\geq0}\|y(t)\|$ (in any $p$-norm).
Through examples we show how the new bounds can be applied. For a large-scale platooning problem, we compute bounds on the deviations from equilibrium for a worst-case bounded initial condition.   

The remainder of this paper is organized as follows. In \secref{Preliminaries} we introduce the preliminaries of this work. Lower and upper bounds on the transient of $\sup_t \|y(t)\|$ and simple examples illustrating how to apply the bounds are presented in \secref{mainresults}.
Then we illustrate an application of our results in the form of vehicle strings in \secref{aplications}. Lastly our conlusions are presented in \secref{conclusions}.

\section{Preliminaries}
\label{sec:Preliminaries}
Consider the linear time-invariant system
\begin{equation}
    \begin{aligned}
        \dot{\xi}(t)&= \mathcal{A}\xi(t) +\mathcal{B}u(t)\\
        y(t)&=\mathcal{C} \xi(t),
    \end{aligned}
    \label{eq:linearsystem}
\end{equation} 
where the state $\xi \in \mathbb{R}^{N}$, $\mathcal{A}\in \mathbb{R}^{N\times N}$, $\mathcal{B}\in \mathbb{R}^{N\times P}$, $\mathcal{C}\in \mathbb{R}^{Q \times N}$, and output $y\in \mathbb{R}^{Q} $. The initial condition is $\xi(0)=\xi_0$. We will interpret $\mathcal{C}(sI-\mathcal{A})^{-1}\mathcal{B}$ as a transfer matrix and call the system~\eqref{eq:linearsystem} input-output stable if all poles of this transfer matrix lie in the open left half plane. 
Denote by $\sigma(\mathcal{A})$ the spectrum, i.e., the set of eigenvalues of $\mathcal{A}$. 



We will often let the system in~\eqref{eq:linearsystem} model 
matrix differential equations of the form
\begin{equation}
    \begin{aligned}
        x^{(l)}(t)+A_{l-1}x^{(l-1)}(t) +\dots +A_{0}x(t) =B u(t)\\
        y(t)=\mathcal{C} \xi(t),
    \end{aligned}
    \label{eq:matpolsys}
\end{equation} 
where $x(t) \in \mathbb{R}^{n}$ and $x^{(k)}$ denotes the $k^{\mathrm{th}}$ time derivative of $x$: $x^{(k)}(t) = \frac{\mathrm{d}^kx(t)}{\mathrm{d}t^k}$. In this case, $\xi(t) = [x, \dot{x},\ldots,x^{(l-1)}]^{\top} \in \mathbb{R}^{nl}$, with $nl=N$.
%
This system can be equivalently stated on block-companion form as
\begin{equation}
\begin{aligned}
    \dot{\xi}(t) &=
    \underbrace{\begin{bmatrix}
        0 &  I_n &0 &\dots  \\ 
        \vdots&\ddots & \ddots &0  \\
        0 &\dots  &0 & I_n\\
        -A_{0}& -A_{1} & \dots & -A_{l-1}
    \end{bmatrix}}_{\mathcal{A}}
    \xi(t)
     +\underbrace{\begin{bmatrix} 0\\ \vdots \\ 0\\ B  \end{bmatrix}}_{\mathcal{B}}u(t)
    \\ 
    y(t)&= \mathcal{C} \xi(t).
    \end{aligned}
    \label{eq:matpol_linear}
\end{equation}


\subsection{Signal and system norms}
Norms are central to this work. Here we will consider the standard vector $p$-norms:
\begin{equation*}
\begin{aligned}
    \|x\|_p&=\begin{cases}
  \left(\sum_{k=1}^N |x_k|^p\right)^{\frac{1}{p}}& \text{ if } 1\leq p<\infty \\
 \max_k |x_k|& \text{ if } p=\infty,
\end{cases}
\end{aligned}
\end{equation*}
where $x\in \mathbb{C}^{N}$. For matrices we consider the corresponding induced norms, i.e.
$$\|\mathcal{A}\|=\sup_{\|x\|=1}\|\mathcal{A} x\|,$$
where $\mathcal{A}\in \mathbb{C}^{M\times N}$.

In general, our results can be interpreted in any of these norms and we will often omit the subscript to indicate that the results are valid for all of them. What we need for our theorems is, more specifically, that the matrix norms are submultiplicative, which means that the following inequality is valid for any two compatible matrices $A_1,~A_2$ 
$$\|A_1 A_2\|\leq \|A_1\|\|A_2\|.$$
It is well known that this is true for all the $p$-norms.

\subsection{Input-output scenarios}
\label{sec:ioscenarios}
We will present bounds in terms of the scaled exponential matrix $\mathcal{C}e^{\mathcal{A}t}\mathcal{B}$. Its norm can be seen as bounds on the transient response of the system~\eqref{eq:linearsystem} in the following scenarios: 
\subsubsection{Impulse response}
Consider the input signal $\{u(t)=\delta(t)u_0\}$ with $u_0\in \mathbb{R}^{P}$ and let $||u_0||=1$ in some norm. The solution of \eqref{eq:linearsystem} is given by 
\begin{equation}
    y(t)=\mathcal{C}e^{t\mathcal{A}}\mathcal{B}u_0
\end{equation}
and the worst possible transient of $y(t)$ is given by $$\sup_{t}\|y(t)\|=\sup_{t}\sup_{\|u_0\|=1}\|\mathcal{C}e^{t\mathcal{A}}\mathcal{B}u_0\|=\sup_{t}\|\mathcal{C}e^{t\mathcal{A}}\mathcal{B}\|.$$

\subsubsection{Response to an initial condition}
An initial condition response is  given by 
$$ y(t)=\mathcal{C}e^{t \mathcal{A}}\xi(0).$$
To study the worst possible initial condition with respect to resulting deviations in the output $y(t)$ we may consider
$$\sup_{t}\|y(t)\|=\sup_{t}\sup_{\|\xi_0\|=1}\|\mathcal{C}e^{t\mathcal{A}}\xi_0\|=\sup_{t}\|\mathcal{C}e^{t\mathcal{A}}\|.$$
The corresponding analysis for the worst-case \textit{structured} initial condition is done by multiplying $\xi_0$ by $\mathcal{B}$. In this case, 
$$\sup_{t}\|y(t)\|=\sup_{t}\sup_{\|\xi_0\|=1}\|\mathcal{C}e^{t\mathcal{A}}\mathcal{B}\xi_0\|=\sup_{t}\|\mathcal{C}e^{t\mathcal{A}}\mathcal{B}\|.$$ For example, $\mathcal{B}=(I,0,\dots,0)^T$ in \eqref{eq:matpol_linear} corresponds to all initial derivatives being zero. 




\subsection{Complex analysis}

The basis for our upcoming theorems is three Laplace transform results, which were also used to derive key results in~{\cite{trefethenbook}}. For completeness they are also presented.

\begin{lemma}[{\cite[Theorem 15.1]{trefethenbook} }]
    \label{thrm:trefethen}
    Let $\mathcal{A}$ be a matrix. There exist $\omega\in \mathbb{R}$ and $M\geq 1$ such that
    \begin{equation}
        \|e^{t \mathcal{A}}\|\leq M e^{\omega t}\qquad \forall t\geq 0.
        \label{eq:expbound}
    \end{equation}
    Any $s\in \mathbb{C}$ with $\text{Re} s>\omega$ is in the resolvent set of $\mathcal{A}$, with 
    \begin{equation}
        (sI-\mathcal{A})^{-1}=\int_0^\infty e^{-s t} e^{t \mathcal{A}} \mathrm{d}t.
        \label{eq:exptoinv}
    \end{equation}
    If $\mathcal{A}$ is a matrix or bounded operator, then
    \begin{equation}
        e^{t \mathcal{A}} =\frac{1}{2 \pi i}\int_\Gamma e^{s t} (sI-\mathcal{A})^{-1} \mathrm{d}s,
        \label{eq:invtoexp}
    \end{equation}
    where $\Gamma$ is any closed and positively oriented contour that encloses $\sigma(\mathcal{A})$ once in its interior.
\end{lemma}

\subsection{Pseudospectra}
Pseudospectra have proven themselves to be a useful tool for analysing the transient behavior and robustness of differential equations, see e.g. \cite{GREEN2006transient_response}. There are several equivalent definitions of the pseudospectra of 
a matrix $\mathcal{A}\in \mathbb{C}^{N\times N}$. Two equivalent and well known are:
\begin{definition}[$\epsilon$-pseudospectra]
    \label{def:eps_spectra}
    \begin{equation}
    \sigma_\epsilon(\mathcal{A})=\{s\in\mathbb{C} \mid \|(sI-\mathcal{A})^{-1}\|>\epsilon^{-1} \}
    \label{eq:pseudo_def_res}
    \end{equation}
\end{definition}
and
\begin{definition}[$\epsilon$-pseudospectra]
    \begin{multline}
    \sigma_\epsilon(\mathcal{A})=\{s\in\mathbb{C} \mid s\in \sigma(\mathcal{A}+E) \ldots \\ \ldots \text{ for some }  E\in \mathbb{C}^{N\times N} \text{ with } \|E\|<\epsilon \},
        \label{eq:pseudo_def_rob}
\end{multline}
\vspace{-5mm}
\end{definition}
where $\sigma(\mathcal{A})$ denotes the (usual) spectrum of a matrix $\mathcal{A}$. We will also make use of the $\epsilon$-pseudospectral abscissa, defined as $\alpha_\epsilon=\sup_{s\in\sigma_\epsilon} \text{Re}s$.

From the two definitions of $\sigma_\epsilon$ we can get an idea of what they are used for. The first relates to the size of the resolvent and enables complex analysis in line with Lemma~\ref{thrm:trefethen}. The latter relates to the robustness of the matrix under perturbations. By considering level curves of pseudospectra for various $\epsilon$-levels it is possible to get an understanding of the solutions of the linear differential equation $\dot{x}(t)=\mathcal{A}x(t)$ and of how sensitive the system is to perturbations.

When one is concerned with the transient behaviour of an input-output system as defined in \eqref{eq:linearsystem} it will be proven useful to generalize Definition \ref{def:eps_spectra} in the following way:
\begin{definition}[Input-output $\epsilon$-pseudospectra]
    \label{def:eps_spectra_new}
    \begin{equation}
    \sigma_\epsilon(\mathcal{A},\mathcal{B},\mathcal{C})=
    \{s\in\mathbb{C} \mid \|\mathcal{C}(sI-\mathcal{A})^{-1}\mathcal{B}\|>\epsilon^{-1} \}.
    \label{eq:io_pseudo_def}
)    \end{equation}
\end{definition}
\vspace{2mm}
The corresponding input-output \emph{pseudospectral abscissa} we define as $$\alpha_\epsilon(\mathcal{A},\mathcal{B},\mathcal{C})=\sup_{s \in \sigma_\epsilon(\mathcal{A},\mathcal{B},\mathcal{C})} \text{Re}(s).$$
We also define the input-output spectrum $\sigma(\mathcal{A},\mathcal{B},\mathcal{C})$ as the set of poles of the transfer matrix $\mathcal{C}(sI-\mathcal{A})^{-1}\mathcal{B}$.



\subsection{Kreiss theorem}
The transient behavior of a matrix exponential for a stable matrix $\mathcal{A}\in \mathbb{R}^{N\times N}$ can be bounded through the so called \textit{Kreiss bounds} \cite[Thrm. 18.5]{trefethenbook}:
\begin{equation}
    \mathcal{K}(\mathcal{A}) \leq \sup_{t\geq 0} \|e^{t \mathcal{A}}\| \leq e N \mathcal{K}(\mathcal{A}).
    \label{eq:Kreiss_standard}
\end{equation}
Here the Kreiss constant is defined as
\begin{equation}
\label{eq:kreissdef}
    \mathcal{K}(\mathcal{A}) =\sup_{\text{Re} s>0}\text{Re} s \|(sI-\mathcal{A})^{-1}\|.
\end{equation}

Comparing~\eqref{eq:kreissdef} to Definition 1, the relation between the Kreiss bound and pseudospectra becomes evident. In fact, it holds that $\mathcal{K}(\mathcal{A}) = \sup_{\epsilon>0}\alpha_\epsilon /\epsilon$.
In a controls context, it is natural to not only consider the matrix exponential, but rather an input-output setting. We therefore define 
the input-output Kreiss constant 
\begin{equation}
    \label{eq:kreissabc}
     \mathcal{K}(\mathcal{A},\mathcal{B},\mathcal{C})=\sup_{\text{Re} s>0}\text{Re} s \|\mathcal{C}(sI-\mathcal{A})^{-1}\mathcal{B}\|=\sup_{\epsilon>0}\frac{\alpha_\epsilon(\mathcal{A},\mathcal{B},\mathcal{C)}}{\epsilon}.
\end{equation}


\section{Input-Output Transient Bounds}
\label{sec:mainresults}
We now make use of the theory in the previous section to derive bounds on the transient performance of the system~\eqref{eq:linearsystem}, under the input-output scenarios introduced earlier. We will give both lower and upper bounds. As a starting point, consider the following proposition, which is a simple but important
extension to Lemma~\ref{thrm:trefethen}: 
\begin{proposition}
\label{prop:boundsintegrals}
    Let $\mathcal{A}$, $\mathcal{B}$ and $\mathcal{C}$ be matrices and let $\|\cdot\|$ denote a submultiplicative norm. There exist $w\in \mathbb{R}$ and $M\geq \|\mathcal{C}\mathcal{B}\|$ such that 
    \begin{equation}
        \|\mathcal{C}e^{t \mathcal{A}}\mathcal{B}\|\leq M e^{\omega t}\qquad \forall t\geq 0.
        \label{eq:CexpB}
    \end{equation}
    Any $s\in \mathbb{C}$ with $\text{Re} s>\omega$ is in the resolvent set of $\mathcal{A}$, with 
    \begin{equation}
        \mathcal{C}(sI-\mathcal{A})^{-1}\mathcal{B}=\int_0^\infty e^{-s t} \mathcal{C}e^{t \mathcal{A}}\mathcal{B} \mathrm{d}t,
        \label{eq:CresBlexp}
    \end{equation}
    \begin{equation}
        \mathcal{C}e^{t \mathcal{A}}\mathcal{B} =\frac{1}{2 \pi i}\int_\Gamma e^{s t}  \mathcal{C}(sI-\mathcal{A})^{-1}\mathcal{B}\mathrm{d}s,
        \label{eq:CexpBlint}
    \end{equation}
    and where $\Gamma$ is any closed and positively oriented contour that encloses $\sigma(\mathcal{A},\mathcal{B},\mathcal{C})$ once in its interior.
\end{proposition}
\begin{proof}
    First, \eqref{eq:CexpB} follows from the norm's submultiplicativity and~\eqref{eq:expbound} as
    \begin{equation*}
        \begin{aligned}
             \|\mathcal{C}e^{t \mathcal{A}}\mathcal{B}\|&\leq \|\mathcal{C}\|\|\mathcal{B}\| \|e^{t \mathcal{A}}\| \leq \|\mathcal{C}\|\|\mathcal{B}\| \hat{M} e^{w t},
         \end{aligned}
    \end{equation*}
    with $M=\|\mathcal{C}\|\|\mathcal{B}\|\hat{M}$. Letting $t=0$ yields $\|\mathcal{C} \mathcal{B}\|\leq M $.
    
  Next,  \eqref{eq:CresBlexp} follows from linearity of the integral, i.e., the fact that for any compatible matrices $B$,  $C$, and $f(x)$ we have $B\int (f(x)\mathrm{d}x)C = \int B f(x)C\mathrm{d}x$.
  
  Last, consider \eqref{eq:CexpBlint}. Through linearity and \eqref{eq:invtoexp}, we get
  \begin{equation*}
      \mathcal{C}e^{t \mathcal{A}}\mathcal{B} =\frac{1}{2 \pi i}\int_{\Gamma'} e^{s t}  \mathcal{C}(sI-\mathcal{A})^{-1}\mathcal{B}\mathrm{d}s,
  \end{equation*}
  where $\Gamma'$ encircles $\sigma(A)$. If $\sigma(\mathcal{A},\mathcal{B},\mathcal{C})=\sigma(\mathcal{A})$ we are done. If not, suppose that there are $n_p$ distinct poles $s_p\in \sigma(A)$ such that $s_p\notin\sigma(\mathcal{A},\mathcal{B},\mathcal{C})$. Let $\Gamma'$ be the union of $\Gamma$ and $n_p$ disjoint circles with radius $\epsilon$ with the poles~$s_p$ at the center. Let $\epsilon$ be sufficiently small such that the $\epsilon$-circles are disjoint from $\sigma(\mathcal{A},\mathcal{B},\mathcal{C})$. Now, since the transfer matrix $\mathcal{C}(sI-\mathcal{A})^{-1}\mathcal{B}$ does not contain any poles in the interior of the $\epsilon$-discs, 
  each of the transfer functions is holomorphic in each disc enclosed by the $\epsilon$-circles. By the maximum modulus principle they cannot have any strict local maximum in the interior of each $\epsilon$-disc. This implies that there is an $M_\epsilon \geq 0$ such that each transfer function $|(\mathcal{C}(sI-\mathcal{A})^{-1}\mathcal{B})_{i,j}|\leq M_\epsilon$.
  In turn, this implies that  $\|e^{st}(\mathcal{C}(sI-\mathcal{A})^{-1}\mathcal{B})\|_\infty\leq e^{Re(s_p+\epsilon) t}M_\epsilon P$ on any circle~$\gamma$, where $P$ is the number of columns of $\mathcal{B}$. The curve integral is thus bounded by $\int_\gamma\|\mathcal{C}(sI-\mathcal{A})^{-1}\mathcal{B}\|_\infty \mathrm{d}s\leq P e^{Re(s_p+\epsilon) t}M_\epsilon \epsilon 2\pi$, which then converges to $0$ as $\epsilon \to 0$. This is true for all $n_p$ circles and so we can ignore the part encircling the non-observable poles. By equivalence of norms, this is true for any $p$-norm. 
\end{proof}

We will now make use of Proposition~\ref{prop:boundsintegrals} to state upper and lower bounds on the quantity $\sup_{t\geq 0}\| \mathcal{C}e^{t \mathcal{A}}\mathcal{B}\|$.


\subsection{Lower bound}
We begin by stating a lower bound analogous to the lower bound in the Kreiss theorem~\eqref{eq:Kreiss_standard}. Despite its relevance to control systems, this extension of the Kreiss theorem has, to our knowledge, not been observed in the literature apart from~\cite{matsuo1994maxamplitude} and \cite{plischke2005transient}. The short proof we present here, however, is new. 
\begin{theorem}[Lower bound]
\label{thrm:lowbound}
\begin{equation}
     \sup_{t\geq0}\| \mathcal{C}e^{t \mathcal{A}}\mathcal{B}\|\geq \sup_{\text{Re} s>0} \text{Re} s\|\mathcal{C}(sI-\mathcal{A})^{-1}\mathcal{B}\|
\end{equation}
\end{theorem}
\vspace{2mm}
\begin{proof}
    Let $M=\sup_{t\geq0} \| \mathcal{C}e^{t \mathcal{A}}\mathcal{B}\|$ and $\text{Re}s>0$. From~\eqref{eq:CresBlexp} we have
    \begin{equation*}
        \begin{aligned}
                \|\mathcal{C}(sI-\mathcal{A})^{-1}\mathcal{B} \|&=\|\int_0^\infty e^{-s t} \mathcal{C}e^{t \mathcal{A}}\mathcal{B} \mathrm{d}t\|\\
        \implies    \|\mathcal{C}(sI-\mathcal{A})^{-1}\mathcal{B}\| &\leq M\int_0^\infty e^{-t\text{Re}s}\mathrm{d}t = \frac{M}{\text{Re}s},
        \end{aligned}
    \end{equation*}
    multiplying both sides by $\text{Re}s$ proves the inequality.
\end{proof}

The theorem reveals that the input-output Kreiss constant~$\mathcal{K}(\mathcal{A},\mathcal{B},\mathcal{C})$ defined in~\eqref{eq:kreissabc} can lower bound  the transient of the system~\eqref{eq:linearsystem} under the input scenarios in \secref{ioscenarios}.

\subsection{Upper bounds}
Now we present three ways to bound the transient from above, again, using Proposition~\ref{prop:boundsintegrals} as a basis for the proofs.


\begin{theorem}[First upper bound]
\label{thrm:upperbound1}
    If $\mathcal{A}$, $\mathcal{B}$, $\mathcal{C}$ are matrices and $L_\epsilon$ is the arc length of the boundary of $\sigma_\epsilon(\mathcal{A},\mathcal{B},\mathcal{C})$
    or of its convex hull for some $\epsilon >0$, then 
    \begin{equation}
        \|\mathcal{C}e^{t \mathcal{A}}\mathcal{B}\| \leq \dfrac{L_\epsilon e^{t\alpha_\epsilon(\mathcal{A},\mathcal{B},\mathcal{C})}}{2 \pi \epsilon},
        \label{eq:upper1new}
    \end{equation}
    where $\alpha_\epsilon(\mathcal{A},\mathcal{B},\mathcal{C}) = \sup \{{\text{Re} s} \mid \|\mathcal{C}(s-\mathcal{A})^{-1}\mathcal{B}\|\leq \epsilon^{-1}\} $.
\end{theorem}
\vspace{1mm}
\begin{proof}
    For any closed contour $\Gamma$ enclosing $\sigma(\mathcal{A},\mathcal{B},\mathcal{C})$ we have \eqref{eq:CexpBlint}.  Taking the norm on both sides gives
            $$\begin{aligned}
                 \|\mathcal{C}e^{t \mathcal{A}}\mathcal{B}\| &= \norm{ \frac{1}{2 \pi i}\int_\Gamma e^{s t} \mathcal{C}(sI-\mathcal{A})^{-1}\mathcal{B} \mathrm{d}s }\\
                 &\leq \frac{1}{2 \pi} \int_\Gamma \|e^{s t} \mathcal{C}(sI-\mathcal{A})^{-1}\mathcal{B} \|\mathrm{d}s
                 \leq \dfrac{L_\epsilon e^{t\alpha_\epsilon(\mathcal{A},\mathcal{B},\mathcal{C})}}{2 \pi \epsilon}.
            \end{aligned}$$
    The second inequality follows since $\|\mathcal{C}(sI-\mathcal{A})^{-1}\mathcal{B} \|\leq \epsilon^{-1}$ along $\Gamma$. The convex hull can be used to reduce the length $L_\epsilon$ of $\Gamma$. This is possible as  $\|\mathcal{C}(sI-\mathcal{A})^{-1}\mathcal{B} \|\leq \epsilon^{-1}$ on the boundary of the convex hull. 
\end{proof}

Theorem~\ref{thrm:upperbound1} is a fairly straightforward extension of~\cite[Theorem 15.2]{trefethenbook}. However, we next present a novel alternative characterization which will prove useful, in particular for classes of higher-order matrix differential equations.  


\begin{theorem}[Second upper bound]
\label{thrm:upper_with_circle}
Let the system~\eqref{eq:linearsystem} with $(\mathcal{A},\mathcal{B},\mathcal{C})$ be input-output stable and let $R=a\|\mathcal{A}\|$ for some $a>1$. Then
    \begin{equation}
    \begin{aligned}
         \|\mathcal{C}e^{t \mathcal{A}}\mathcal{B}\| &\leq \frac{1}{2 \pi} \!\! \int_{-R}^R \! \!\| \mathcal{C}(i\omega -\mathcal{A})^{-1}\mathcal{B} \|d\omega +\frac{\|\mathcal{C}\| \| \mathcal{B}\|}{2-2a^{-1}}\\
     \end{aligned}
     \label{eq:upper_with_circle}
    \end{equation}
\end{theorem}
\vspace{2mm}
\begin{proof}
By definition of input-output stability all poles of the transfer matrix $\mathcal{C}(sI-\mathcal{A})^{-1}\mathcal{B}$ lie in the left half plane. Furthermore, the spectrum $\sigma(\mathcal{A})$ is contained in the disc $|s|\leq \|\mathcal{A}\|$ since for any eigenvector $x$ of $\mathcal{A}$ with $\|x\|=1$ we have $\|\mathcal{A}\|\geq \|\mathcal{A}x\|=|\lambda|$. Now take $\Gamma$ to be the semicircle with radius $R>\|A\|$ that goes up the imaginary axis and then extends into the left half plane. Then this $\Gamma$ encloses the input-output spectrum $\sigma(\mathcal{A},\mathcal{B},\mathcal{C})$ and \eqref{eq:CexpBlint} yields 
\begin{equation*}
    \begin{aligned}
         \|\mathcal{C}e^{t \mathcal{A}}\mathcal{B}\|&\leq \frac{1}{2 \pi } \int_{\Gamma} e^{t \text{Re}(s)} \| \mathcal{C}(s -\mathcal{A})^{-1}\mathcal{B} \|\mathrm{d}s \\
         &\leq \frac{1}{2 \pi } \int_{-R}^R \| \mathcal{C}(i\omega -\mathcal{A})^{-1}\mathcal{B}\| \mathrm{d}\omega \\
         & + \frac{1}{2 \pi }\int_{\pi/2}^{3\pi/2} \| \mathcal{C}(R e^{i\theta} -\mathcal{A})^{-1}\mathcal{B}\| R \mathrm{d}\theta.\\
     \end{aligned}
\end{equation*}
If $|s|=\|\mathcal{A}\|a$ and $a>1$, then, the integral in the second term can be bounded using the following series expansion of the inverse:
\begin{equation*}
    \|(sI-\mathcal{A})^{-1}\| =\norm{\frac{1}{s} \sum_{k=0}^\infty \left(\frac{\mathcal{A}}{s}\right)^k } 
\leq \frac{1}{\|\mathcal{A}\|}\frac{1}{a-1}
\end{equation*}
This, together with submultiplicativity yield $\| \mathcal{C}(R e^{i\theta} -\mathcal{A})^{-1}\mathcal{B}\| R\leq \|\mathcal{C}\|\mathcal{B}\|\|\mathcal{A}\|a \frac{1}{\|\mathcal{A}\|}\frac{1}{a-1}$ which can be used to upper bound the second integral to $\frac{\|\mathcal{C}\| \| \mathcal{B}\|}{2-2a^{-1}}$.
\end{proof}
\vspace{0.5mm}

Now we will look into another bound, similar in its nature.
\begin{theorem}[Third Upper bound]
\label{thrm:upper_im_axis}
Let the system~\eqref{eq:linearsystem} with $(\mathcal{A},\mathcal{B},\mathcal{C})$ be input-output stable.  If $\|\mathcal{C}(sI-\mathcal{A})^{-1}\mathcal{B}\|\leq M|s|^{-\beta}$ for all $|s|\geq K$ for some $\beta>1$, $M>0$, and $K>0$. Then 
    \begin{equation}
        \|\mathcal{C}e^{t \mathcal{A}}\mathcal{B}\| \leq \frac{1}{2 \pi}\int_{-\infty}^\infty \|\mathcal{C}(i\omega-\mathcal{A})^{-1}\mathcal{B} \| \mathrm{d}\omega<\infty.
        \label{eq:upper2new}
    \end{equation}
\end{theorem}
\vspace{2mm}

\begin{proof}
Taking the norm of~\eqref{eq:CexpBlint}, we get
    \[ \|\mathcal{C}e^{t \mathcal{A}}\mathcal{B}\| \leq \frac{1}{2 \pi}\int_\Gamma \|e^{s t} \mathcal{C}(sI-\mathcal{A})^{-1}\mathcal{B} \| \mathrm{d}s.\]
    Now, use the same semicircle $\Gamma$ as in the proof of Theorem~\ref{thrm:upper_with_circle} with radius $R>\|\mathcal{A}\|$. This $\Gamma$ encloses the input-output spectra $\sigma(\mathcal{A},\mathcal{B},\mathcal{C})$. Furthermore if $R>K$ we have
    \begin{multline*}
         \|\mathcal{C}e^{t \mathcal{A}}\mathcal{B}\| \leq \\ \frac{1}{2 \pi} \lim_{R\rightarrow \infty} \left(  \int_{-R}^R \|e^{i\omega t} \mathcal{C}(i\omega -\mathcal{A})^{-1}\mathcal{B} \|d\omega +\pi R M R^{-\beta} \right)\\
         =\frac{1}{2 \pi} \int_{-\infty}^\infty \| \mathcal{C}(i\omega -\mathcal{A})^{-1}\mathcal{B} \|d\omega,
    \end{multline*}
    where the last equality follows from the condition $\beta>1$. 
\end{proof}

The condition on $\beta$ in Theorem \ref{thrm:upper_im_axis} can be related to the relative degree of the system. For instance, if $\mathcal{C}(sI-\mathcal{A})^{-1}\mathcal{B}=(s^2I+sA_1+A_0)^{-1}$ and is input-output stable, then $\beta=2$ and it is possible to apply the theorem.


The usefulness 
of the three upper bounds boils down to the 
fact that the spectrum is usually difficult to characterize. 
For the first bound~\eqref{eq:upper1new}, 
a good description of the pseudospectra is needed, while in the second and third bounds~\eqref{eq:upper_with_circle}--\eqref{eq:upper2new} good knowledge of the resolvent along the imaginary axis is needed.
We will clarify through two simple examples.

\begin{example}
    Consider the dynamical system
    \begin{equation}
    \begin{aligned}
    \dot{\xi}&=\begin{bmatrix}
        0 & 1\\
        -1 & -2
        \end{bmatrix}
        \xi+\begin{bmatrix}
        0 \\
        1
        \end{bmatrix}u\\
        y&=\begin{bmatrix}
        1 & 0
        \end{bmatrix}\xi.
    \end{aligned}
    \label{eq:simpleexamplesystem}
        \end{equation}
Suppose we are interested in the impulse response of the system. Then we have
\begin{equation*}
    \mathcal{C}(sI-\mathcal{A})^{-1}\mathcal{B}=\frac{1}{(s+1)^2}.
\end{equation*}
In this case we can see that the $\epsilon$-level curves of $\|\mathcal{C}(sI-\mathcal{A})^{-1}\mathcal{B}\|=1/\epsilon$ are given by the circles $|s+1|=\sqrt{\epsilon}$. From~\eqref{eq:upper1new} we see that the upper bound for each $\epsilon$ is
$$ \|\mathcal{C} e^{tA}\mathcal{B}\|\leq \dfrac{2\pi \sqrt{\epsilon} e^{t(-1+\sqrt{\epsilon})}}{2 \pi \epsilon}=\dfrac{e^{t(-1+\sqrt{\epsilon})}}{\sqrt{\epsilon}}.$$ The lowest upper bound is achieved for $\epsilon=1$ and is simply $\|C e^{t \mathcal{A}}\mathcal{B}\| \leq 1$.

The third upper bound (Theorem~\ref{thrm:upper_im_axis}) requires input-output stability, which is clearly satisfied. The relative degree is $2$ which implies $\beta=2>1$. To calculate the bound~\eqref{eq:upper2new} we need to calculate the integral along the imaginary axis. In this case
\begin{align*}
    \|\mathcal{C} e^{t \mathcal{A}}\mathcal{B}\| &\leq\frac{1}{2\pi} \int_{-\infty}^\infty \|\mathcal{C}(\omega i I-\mathcal{A})^{-1}\mathcal{B}\|\mathrm{d}\omega\\
    &=\frac{1}{\pi} \int_{0}^\infty \frac{1}{\omega^2+1}\mathrm{d}\omega =\frac{1}{2}
\end{align*}
A lower bound of this system can be calculated by only considering the real axis (and in this case this is also optimal). This leads to optimizing 
\[\sup_{t\geq 0}\|\mathcal{C} e^{t \mathcal{A}}\mathcal{B}\| \geq \sup_{x>0}\norm{\frac{x}{(x+1)^2}}=\frac{1}{4}.\]
Since this system is very simple it is also possible to calculate the actual maximum which is
$\sup_{t\geq 0}\|\mathcal{C} e^{t \mathcal{A}}\mathcal{B}\|= 1/e.$ 
\end{example}
\vspace{1.2mm}

As demonstrated above, Theorem \ref{thrm:upper_im_axis} is useful if the relative degree of the system transfer function is greater than~1. Now we show a case where where we cannot use this theorem.


\begin{example}
    Consider the same system~\eqref{eq:simpleexamplesystem} as before, but now the response to a non-zero initial value $\xi(0) = [x_0,0]^\top$. This can be represented by $\mathcal{B}=[1,0]^T$. Then we have
    \begin{equation*}
        \mathcal{C}(sI-\mathcal{A})^{-1}\mathcal{B}=\frac{s+2}{(s+1)^2}.
    \end{equation*}
    To calculate our first upper bound in \eqref{eq:upper1new} we need to encircle the spectrum. The shape here is non-trivial but at least we know that for any circle centered around $-1$ with radius smaller than $1$ we have $$\left|\frac{s+2}{(s+1)^2}\right|\leq \left|\frac{2}{(s+1)^2}\right|.$$
    The previous calculations give us the upper bound $\|\mathcal{C} e^{t \mathcal{A}}\mathcal{B}\|\leq 2$. 
    In this case we cannot use Theorem~\ref{thrm:upper_im_axis} since the relative degree is $1$ and therefore $\beta\leq 1$. However, we can use the very similar Theorem~\ref{thrm:upper_with_circle} to give an upper bound. $\|\mathcal{A}\|_\infty=3$ and so for any $R>3$ we can use the theorem. It remains to calculate the curve integral along the imaginary axis. Doing this with numerical integration for $R=9$ yields the upper bound $\|\mathcal{C} e^{tA}\mathcal{B}\|\lessapprox 2.025$ (Through optimization this bound can be lowered to $\|\mathcal{C} e^{tA}\mathcal{B}\|\lessapprox2.023$)
            \end{example}
\vspace{1.5mm}

    Through these examples we have shown that the best upper bound depends on the situation. The upside of using Theorems~\ref{thrm:upper_with_circle} and~\ref{thrm:upper_im_axis} is that they are quite easy to compute numerically. Theorems~\ref{thrm:lowbound} and~\ref{thrm:upperbound1} relate to the level curves of the input-output pseudospectra and can be qualitatively seen through inspection of these curves, as we will demonstrate in the next section. 

\section{Application to networks: vehicle strings}
\label{sec:aplications}
To illustrate our bounds, we consider the problem of controlling a string of vehicles -- the platooning problem. While performance bounds on platoons and their relation to the network or interaction structure has received ample attention, as we stated in the introduction, the problem calls for
bounds relating to the quantity $\sup_{t} \|y(t)\|_\infty$, where $y$ captures a displacement error. 
Bounds of this type are important, especially in platooning, since they directly relate to the allowable spacing between consecutive vehicles. However, they tend to be difficult to derive analytically.  
Here we illustrate how our pseudospectra-inspired approach can be used to evaluate string stability properties for various platoon structures in terms of this quantity. 

For this purpose, consider a platoon of size $n$ where each unit is modelled as a double integrator in one spatial dimension, i.e.
$$\ddot{x}_k=u_k$$
where $x_k$ is the position of the $k$th vehicle with respect to a {fix reference} 
and $u_k$ is the input force at vehicle $k$.

To control the platoon we consider a control law that depends on relative distances to neighboring vehicles and relative to a speed reference. For $k\in \{2,\dots,n-1\}$ we get:
\begin{multline}
\label{eq:vehcontrol}
    u_k = (1+\beta_d)(\dot{x}_{k-1}-\dot{x}_k) - (1-\beta_d)(\dot{x}_{k}-\dot{x}_{k+1})+ \alpha(v_\text{ref} - 
\dot{x}_k)\\
\!+\!(1\!+\!\beta_p)(x_{k-1}\!-\!x_{k}\!-\!d)\! -\!(1-\beta_p)(x_{k}\!-\!x_{k+1}\!-\!d),
\end{multline}
where $\beta_d$ and $\beta_p$ are parameters capturing the degree of symmetry in the control law (i.e., look-ahead vs. look-behind control), $d$ a desired intervehicle spacing, $v_\text{ref}$ is a velocity reference, and $\alpha\ge 0$ is a weight. 
For the first and last vehicles, we simply define $\ddot{x}_1 = - (1-\beta_d)(\dot{x}_{1}-\dot{x}_{2})+ \alpha(v_\text{ref}-\dot{x}_1) -(1-\beta_p)(x_{1}-x_{2}-d),$ $\ddot{x}_n = (1+\beta_d)(\dot{x}_{n-1}-\dot{x}_{n})+ \alpha(v_\text{ref}-\dot{x}_n) +(1+\beta_p)(x_{n-1}-x_{n}-d).$
By considering the translated dynamics $\hat{x}_k=x_k+ k d$ we get the same dynamics as if we assume $d=0$, so for simplicity we set $d=0$ and consider the dynamics around this equilibrium.

The closed-loop system can be written: 
\begin{align}
\nonumber
    \begin{bmatrix} \dot{x}\\\ddot{x}
    \end{bmatrix} &= 
    \begin{bmatrix} 0 & I \\
                    -L_p & -L_d-\alpha I
    \end{bmatrix}
     \begin{bmatrix} x\\\dot{x}
     \end{bmatrix}
     +\begin{bmatrix} 0\\ \alpha \mathbf{1}
     \end{bmatrix}
     v_{\text{ref}}
     =\mathcal{\mathcal{A}}\xi + \mathcal{B}v_{\text{ref}}.\\
     y&=\begin{bmatrix}C & 0 \end{bmatrix}\begin{bmatrix} x\\ \dot{x}
     \end{bmatrix} =\mathcal{C}\xi,
    \label{eq:global_consensus}
\end{align}
where $L_p,~L_d\in \mathbb{R}^{n\times n}$ are graph Laplacians capturing the vehicle interactions (see further down for definitions). 
This and similar systems are well studied, see e.g.~\cite{studli2017StringConcepts}.


A way to ensure the platoon is well-behaved (e.g. string stable) is to make $\alpha$ large in comparison to $L_d$ and $L_p$. This, however, essentially transforms the problem to an open-loop system, which is obviously problematic in a real-world setting with disturbances and measurement noise or bias. This motivates the use of a fairly small $\alpha$, allowing the inter-vehicle adjustments to dominate. In this example, we will use~$\alpha = 0.1$.

We use the framework from Section~\ref{sec:mainresults} to analyze this system for two cases, one where both Laplacians are asymmetric (a directed string) and one where both are symmetric (bidirectional string). For each system we consider the output 
$$y =\begin{bmatrix}
x_1-x_2\\ x_{\lfloor N/2\rfloor}-x_{\lfloor N/2\rfloor+1}\\ x_{n-1}-x_{n}
\end{bmatrix},
$$
which samples three inter-vehicle distances: at the start, middle, and end of the platoon. We will consider the initial condition response, i.e. $\mathcal{B}=I_{2n}$. We expect a string unstable system to perform poorly for at least one of these outputs. 

One merit of our proposed method is the possibility to analyze very large systems. In both cases considered here, we therefore model a platoon of $n = 400$ vehicles. We remark that it would of course be possible to simulate the systems for many different inputs and through simulation bound the possible outputs. But as $n$ grows, this quickly becomes very computationally heavy. Using our theorems generates bounds on the worst case input without any additional effort.  



\subsection{Directed vehicle string} 
\begin{figure}
    \centering
    \def\svgwidth{0.9\linewidth}
    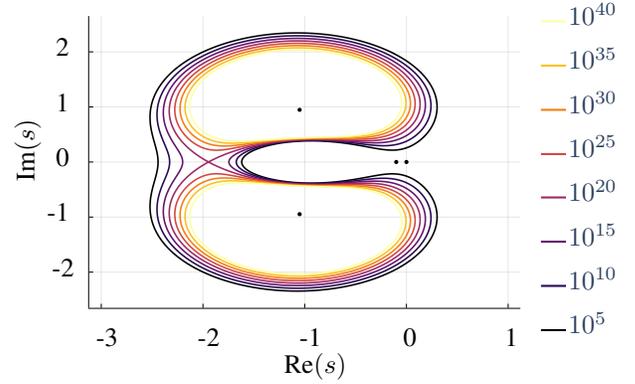
    \caption{The input-output pseudospectra of \eqref{eq:global_consensus} for a directed vehicle string with $n = 400$ vehicles. The black dots are the eigenvalues of $\mathcal{A}$. The large values of the input-output pseudospectra even for small $s$ in the right half plane indicate an unfavorable lower bound in Theorem~\ref{thrm:lowbound}. }
    \label{fig:asym_spectra}
\end{figure}

Consider the control law~\eqref{eq:vehcontrol} with $\beta_p=\beta_d=1$, which renders it fully asymmetric. In this case, we obtain in~\eqref{eq:global_consensus} 
$L_p=L_d=L_\text{asym}$, with
\begin{equation}
    L_\text{asym}=\begin{bmatrix} 0 & & \\
-2 & 2 &  \\
& \ddots & \ddots & \\
& &  -2 & 2 
\end{bmatrix}
.
\label{eq:asymL}
\end{equation}
In \figref{asym_spectra} we show the shape of the input-output pseudospectra corresponding to an initial condition response, that is, the level curves of the quantity $\|\mathcal{C}(sI-\mathcal{A})^{-1}I_{2n}\|_\infty$ of~\eqref{eq:global_consensus} with $n =400$ vehicles. We can see that the level curves extend far into the right half plane with magnitudes of order $10^{30}$ where $\text{Re}(s)$ is of order $10^{-1}$. Through the lower bound in Theorem~\ref{thrm:lowbound} we can immediately see that there will be a large transient of $\|y(t)\|_\infty$ in at least the orders of $10^{29}$. Through a line search, starting at the maximum along the imaginary axis and going into the right half plane we learn that the lower bound in amplification from the worst-case initial conditions to the output (see Section~\ref{sec:ioscenarios}) is at least 
$$\sup_{t\geq 0} \| y(t)\|_\infty \gtrapprox 4.3\cdot 10^{31}, $$ 
for some $\xi_0$ such that $\|\xi_0\|_\infty \leq 1$.

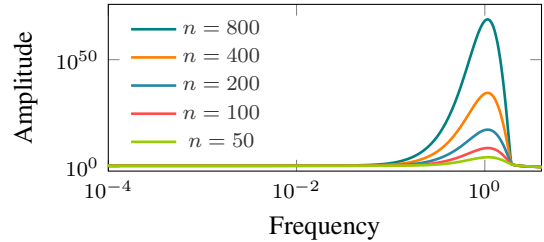
\begin{figure}
    \centering
    \def\svgwidth{2\linewidth}
    \begin{tikzpicture}[trim axis left,trim axis right,baseline]
\begin{axis}[
width=.85\columnwidth,
height=3.8cm,
xlabel={Frequency},
ylabel={Amplitude},
xmode=log,
ymode=log,
xmin=10^-4,
 xmax=4,
ymin=10^-2,
xtick={1,10^-2,10^-4,10^-6},
cycle list name=exotic,
ytick={1,10^50, 10^100},
legend style={fill=white,font=\footnotesize, fill opacity=.8, draw opacity=0.,style={draw=none}},
legend pos=north west,
tick label style={font=\footnotesize},
extra y tick style={grid=major}
]
    
        \addplot +[mark=none,solid,smooth,opacity=1,line width=1pt] table[x index=0, y index=1,col sep=comma] {N800bode_asym.csv};
    \addplot +[mark=none,solid,smooth,opacity=1,line width=1pt] table[x index=0, y index=1,col sep=comma] {N400bode_asym.csv};
    \addplot +[mark=none,solid,smooth,opacity=1,line width=1pt] table[x index=0, y index=1,col sep=comma] {N200bode_asym.csv};
    \addplot +[mark=none,solid,smooth,opacity=1,line width=1pt] table[x index=0, y index=1,col sep=comma] {N100bode_asym.csv};
    \addplot +[mark=none,solid,smooth,opacity=1,line width=1pt] table[x index=0, y index=1,col sep=comma] {N50bode_asym.csv};
    
        \addlegendentry{$n=800$}
    \addlegendentry{$n=400$}
    \addlegendentry{$n=200$}
    \addlegendentry{$n=100$}
    \addlegendentry{$n=50$}


\end{axis}
\end{tikzpicture}
    \caption{Bode plot displaying the worst-case frequency response of~\eqref{eq:global_consensus} for a directed vehicle string. The amplitude is measured in $\|\cdot\|_\infty$ and the response is shown for various platoon lengths $n$.}
    \label{fig:asym_bode}
\end{figure}

\figref{asym_bode} displays a Bode plot of the system for various platoon sizes $n$. That is, we plot the amplitude $\|\mathcal{C}(sI-\mathcal{A})^{-1}I_{2n}\|_\infty$ for $s = i\omega$, $\omega \in (0,\infty)$. Here,
we can see the extreme amplification of the frequency response close to the frequency $\omega =1$. According to our Theorem~\ref{thrm:upper_with_circle}, we can use this frequency response to 
calculate an upper bound on the transient through integration. 
By numerical integration we can estimate the  upper bound to be
$$\sup_{t\geq 0}|y(t)\|_\infty \lessapprox 1.4\cdot 10^{33}. $$ 

From these two bounds we can already conclude that this topology is not suitable for a string of vehicles. 

\begin{rem}
The bounds were possible to compute, since the inversion $\mathcal{C}(sI-\mathcal{A})^{-1}$ can be reduced to the sparse problem $C(s^2I+sL_d+L_p)^{-1}$, where $C\in\mathbb{R}^{3,n}$. As $L_d$ and $L_p$ are tridiagonal this can be computed in $\mathcal{O}(n)$ operations.

\end{rem}


\subsection{Bidirectional (symmetric) vehicle string} 
\begin{figure}
    \centering
    \def\svgwidth{0.9\linewidth}
    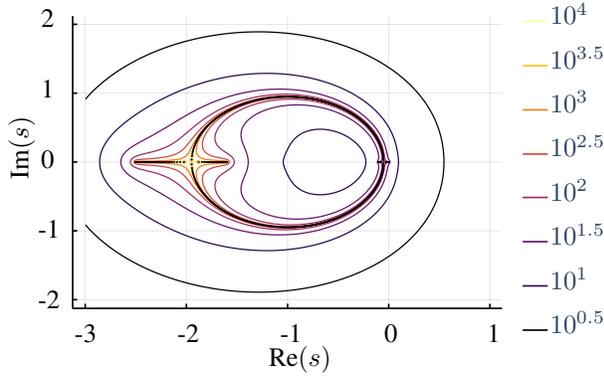
    \caption{The input-output pseudospectra of \eqref{eq:global_consensus} for a bidirectional vehicle string with $n = 400$ vehicles. The black dots are the eigenvalues of $\mathcal{A}$. The level curve corresponding to $\|\mathcal{C}(sI-\mathcal{A})^{-1}\|=10^{1.5}$ can be roughly inscribed in a $1$ radius circle, which hints through  Theorem~\ref{thrm:upperbound1} that the transients of $\|y(t)\|_\infty$ will be small. }
    \label{fig:sym_spectra}
\end{figure}

Now, let $\beta_d=\beta_p = 0$ in~\eqref{eq:vehcontrol}, leading to the symmetric Laplacians $L_p=L_d=L_\text{sym}$, with
    \begin{equation}
        L_\text{sym}=\begin{bmatrix} 1 &-1 & &\\
-1 & 2 & -1 & \\
& \ddots & \ddots & \ddots &\\
& & -1 & 2 & -1\\
& &  & -1 & 1 
\end{bmatrix}.
\label{eq:symL}
    \end{equation}
This corresponds to a bidirectional string of vehicles. In \figref{sym_spectra} we show the shape of the input-output pseudospectra corresponding to an initial condition response of \eqref{eq:global_consensus} with $n=400$ vehicles. We can see that the input-output spectra is quite well-behaved and do not extend far into the right half plane, indicating transients will be modest. 
By Theorem~\ref{thrm:lowbound} and a line search along the real axis, we learn that the lower bound in amplification from initial conditions to output is
$$\sup_{t\geq 0} \| y(t)\|_\infty \gtrapprox 2.2, $$
for some $\xi_0$ such that $\|\xi_0\|_\infty \leq 1$.

\begin{figure}
    \centering
    \def\svgwidth{0.4\linewidth}
    \begin{tikzpicture}[trim axis left,trim axis right,baseline]
\begin{axis}[
width=.85\columnwidth,
height=3.95cm,
xlabel={Frequency},
ylabel={Amplitude},
xmode=log,
ymode=log,
xmin=10^-12,
 xmax=4,
ymax=10^9,
xtick={5,10^-3,10^-6,10^-9, 10^-12},
extra y tick style={grid=major},
cycle list name=exotic,
ytick={1,10^3, 10^6},
legend pos= north east,
tick label style={font=\footnotesize},
legend style={fill=white, font=\footnotesize,fill opacity=.8, draw opacity=0.,style={draw=none}}
]
    \addplot +[mark=none,solid,smooth,opacity=1,line width=1pt] table[x index=0, y index=1,col sep=comma] {N1000000bode_sym.csv};
    \addplot +[mark=none,solid,smooth,opacity=1,line width=1pt] table[x index=0, y index=1,col sep=comma] {N100000bode_sym.csv};
    \addplot +[mark=none,solid,smooth,opacity=1,line width=1pt] table[x index=0, y index=1,col sep=comma] {N10000bode_sym.csv};
    \addplot +[mark=none,solid,smooth,opacity=1,line width=1pt] table[x index=0, y index=1,col sep=comma] {N1000bode_sym.csv};
    \addplot +[mark=none,solid,smooth,opacity=1,line width=1pt] table[x index=0, y index=1,col sep=comma] {N100bode_sym.csv};

    \addlegendentry{$n=10^6$}
    \addlegendentry{$n=10^5$}
    \addlegendentry{$n=10^4$}
    \addlegendentry{$n=10^3$}
    \addlegendentry{$n=10^2$}


\end{axis}
\end{tikzpicture}
    \caption{Bode plot displaying the worst-case frequency response of~\eqref{eq:global_consensus} for a bidirectional vehicle string. The amplitude is measured in $\|\cdot\|_\infty$ and the response is shown for various platoon lengths $n$.}
    \label{fig:sym_bode}
\end{figure}
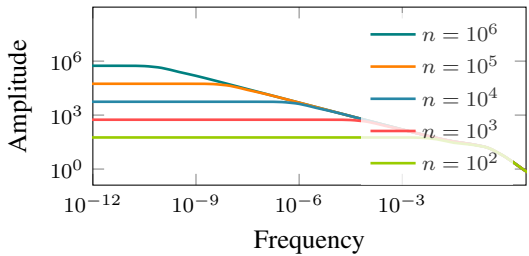

In \figref{sym_bode} we can see the frequency response calculated for various $n$. Interestingly, the common slope among the curves seems to only behave like a square root which would mean that there is an upper bound independent of the platoon length which bounds the transients due to arbitrary non-zero initial conditions. The numerical upper bound when $n=400$ was calculated to
$$\sup_{t\geq0}|\mathcal{C}e^{t\mathcal{A}}\xi_0\|_\infty \lessapprox 9.3. $$ 
This can be compared with the upper bound calculated for $n=10^6$ which evaluated at $9.4$.


\section{Conclusions}
\label{sec:conclusions}

In this work we have proposed a pseudospectra-based approach to analyze transient performance of input-output systems, and generalized existing bounds for this purpose. 
Through our bounds it is possible to quantify $\sup_{t\geq0}\|\mathcal{C}e^{t\mathcal{A}}\mathcal{B}\|$ in terms of lower and upper bounds. These can be seen as bounding the performance from a worst-case input disturbance to an output $y(t)$ in any p-norm -- otherwise often intractable to study. Regarding the problem of controlling vehicle strings in~\secref{aplications}, we illustrated one application where we believe our bounds can be useful, opening the door to future analysis. For instance, deriving analytical bounds for special network structures. The theorems can also be used to numerically calculate bounds for network structures where the worst inputs are non-obvious, for instance when the agents are non-homogeneous or interaction matrices non-normal.

\bibliographystyle{IEEETran}
\bibliography{references}

\begin{thebibliography}{10}
\providecommand{\url}[1]{#1}
\csname url@rmstyle\endcsname
\providecommand{\newblock}{\relax}
\providecommand{\bibinfo}[2]{#2}
\providecommand\BIBentrySTDinterwordspacing{\spaceskip=0pt\relax}
\providecommand\BIBentryALTinterwordstretchfactor{4}
\providecommand\BIBentryALTinterwordspacing{\spaceskip=\fontdimen2\font plus
\BIBentryALTinterwordstretchfactor\fontdimen3\font minus
  \fontdimen4\font\relax}
\providecommand\BIBforeignlanguage[2]{{%
\expandafter\ifx\csname l@#1\endcsname\relax
\typeout{** WARNING: IEEEtran.bst: No hyphenation pattern has been}%
\typeout{** loaded for the language `#1'. Using the pattern for}%
\typeout{** the default language instead.}%
\else
\language=\csname l@#1\endcsname
\fi
#2}}

\bibitem{FaxMurray}
J.~A. Fax and R.~M. Murray, ``Information flow and cooperative control of
  vehicle formations,'' \emph{IEEE Trans. Autom. Control}, vol.~49, no.~9, pp.
  1465--1476, Sept 2004.

\bibitem{OlfatiSaber2004}
R.~Olfati-Saber and R.~M. Murray, ``Consensus problems in networks of agents
  with switching topology and time-delays,'' \emph{IEEE Trans. Autom. Control},
  vol.~49, no.~9, pp. 1520--1533, Sept 2004.

\bibitem{Bamieh2012}
B.~Bamieh, M.~R. Jovanovi\'c, P.~Mitra, and S.~Patterson, ``Coherence in
  large-scale networks: {D}imension-dependent limitations of local feedback,''
  \emph{IEEE Trans. Autom. Control}, vol.~57, no.~9, pp. 2235 --2249, Sept.
  2012.

\bibitem{SiamiMotee2015}
M.~Siami and N.~Motee, ``Fundamental limits and tradeoffs on disturbance
  propagation in large-scale dynamical networks,'' \emph{IEEE Trans. Autom.
  Control}, vol.~61, no.~12, pp. 4055--4062, 2016.

\bibitem{Tegling2019TAC}
E.~Tegling, P.~Mitra, H.~Sandberg, and B.~Bamieh, ``On fundamental limitations
  of dynamic feedback control in regular large-scale networks,'' \emph{IEEE
  Trans. Autom. Control}, vol.~64, no.~12, pp. 4936--4951, 2019.

\bibitem{LevineAthans1966}
W.~Levine and M.~Athans, ``On the optimal error regulation of a string of
  moving vehicles,'' \emph{IEEE Trans. Autom. Control}, vol.~11, no.~3, pp.
  355--361, 1966.

\bibitem{chu1974decentralized_control}
K.~Chu, ``Decentralized control of high-speed vehicular strings,''
  \emph{Transportation Science}, vol.~8, no.~4, pp. 361--384, 1974.

\bibitem{swaroop1996stringstability}
D.~Swaroop and J.~Hedrick, ``String stability of interconnected systems,''
  \emph{IEEE Trans. Autom. Control}, vol.~41, no.~3, pp. 349--357, March 1996.

\bibitem{seiler2004disturbancep_propagation}
P.~Seiler, A.~Pant, and K.~Hedrick, ``Disturbance propagation in vehicle
  strings,'' \emph{IEEE Trans. Autom. Control}, vol.~49, no.~10, pp.
  1835--1842, Oct 2004.

\bibitem{studli2017StringConcepts}
S.~Stüdli, M.~Seron, and R.~Middleton, ``From vehicular platoons to general
  networked systems: String stability and related concepts,'' \emph{Annu. Rev.
  Control}, vol.~44, pp. 157--172, 2017.

\bibitem{FENG2019StringDefs}
S.~Feng, Y.~Zhang, S.~E. Li, Z.~Cao, H.~X. Liu, and L.~Li, ``String stability
  for vehicular platoon control: Definitions and analysis methods,'' \emph{Ann.
  Rev. Control}, vol.~47, pp. 81--97, 2019.

\bibitem{HERMAN2017diffasym}
I.~Herman, S.~Knorn, and A.~Ahlén, ``Disturbance scaling in bidirectional
  vehicle platoons with different asymmetry in position and velocity
  coupling,'' \emph{Automatica}, vol.~82, pp. 13--20, 2017.

\bibitem{ploeg2014Lpstability}
J.~Ploeg, N.~van~de Wouw, and H.~Nijmeijer, ``Lp string stability of cascaded
  systems: Application to vehicle platooning,'' \emph{IEEE Trans. Control Syst.
  Technol.}, vol.~22, no.~2, pp. 786--793, March 2014.

\bibitem{feintuch2012}
A.~Feintuch and B.~Francis, ``Infinite chains of kinematic points,''
  \emph{Automatica}, vol.~48, no.~5, pp. 901--908, 2012.

\bibitem{trefethenbook}
L.~N. Trefethen and M.~Embree, \emph{Spectra and Pseudospectra: The Behavior of
  Nonnormal Matrices and Operators}.\hskip 1em plus 0.5em minus 0.4em\relax
  Princeton: Princeton University Press, 2005.

\bibitem{Kreiss1962}
H.-O. Kreiss, ``{\"U}ber die stabilit{\"a}tsdefinition f{\"u}r
  differenzengleichungen die partielle differentialgleichungen approximieren,''
  \emph{BIT Numerical Mathematics}, vol.~2, pp. 153--181, 1962.

\bibitem{Structured2001Tisseur}
F.~Tisseur and N.~J. Higham, ``Structured pseudospectra for polynomial
  eigenvalue problems, with applications,'' \emph{SIAM J. Matrix Anal. Appl.},
  vol.~23, no.~1, pp. 187--208, 2001.

\bibitem{lancaster2005matrix_polynomials}
P.~Lancaster and P.~Psarrakos, ``On the pseudospectra of matrix polynomials,''
  \emph{SIAM J. Matrix Anal. Appl.}, vol.~27, no.~1, pp. 115--129, 2005.

\bibitem{GREEN2006transient_response}
K.~Green, A.~Champneys, and M.~Friswell, ``Analysis of the transient response
  of an automatic dynamic balancer for eccentric rotors,'' \emph{Int. J. Mech.
  Sci.}, vol.~48, no.~3, pp. 274--293, 2006.

\bibitem{matsuo1994maxamplitude}
T.~Matsuo, ``Some transfer-function conditions for a desired maximum amplitude
  or exponential envelope of a closed-loop transient response,'' in \emph{IEEE
  Conf. on Decision and Control}, vol.~3, 1994, pp. 2659--2660.

\bibitem{plischke2005transient}
E.~Plischke, ``Transient effects of linear dynamical systems,'' Ph.D.
  dissertation, University of Bremen, 2005.

\end{thebibliography}
\vspace{1mm}
\end{document}